\documentclass[11pt]{article}
\usepackage{amsmath,amsthm,amsfonts}
\usepackage{color}
\numberwithin{equation}{section} \allowdisplaybreaks
 \usepackage[all]{xy}

\newtheorem{theorem}{\sc Theorem}[section]
\newtheorem{lemma}[theorem]{\sc Lemma}
\newtheorem{proposition}[theorem]{\sc Proposition}
\newtheorem{corollary}[theorem]{\sc Corollary}
\newtheorem{definition}[theorem]{\sc Definition}
\newtheorem{example}[theorem]{\sc Example}
\newtheorem{examples}[theorem]{\sc Examples}
\newtheorem{remark}[theorem]{\sc Remark}

\newcommand{\bet}{\begin{theorem}}
\newcommand{\eet}{\end{theorem}}
\newcommand{\blm}{\begin{lemma}}
\newcommand{\elm}{\end{lemma}}
\newcommand{\bprop}{\begin{proposition}}
\newcommand{\eprop}{\end{proposition}}
\newcommand{\bcor}{\begin{corollary}}
\newcommand{\ecor}{\end{corollary}}
\newcommand{\bdf}{\begin{definition}\rm}
\newcommand{\edf}{\end{definition}}
\newcommand{\bp}{\begin{proof}}
\newcommand{\ep}{\end{proof}}
\newcommand{\bex}{\begin{example}\rm}
\newcommand{\eex}{\end{example}}
\newcommand{\bremark}{\begin{remark}\rm}
\newcommand{\eremark}{\end{remark}}
\newcommand{\bexs}{\begin{examples}\rm}
\newcommand{\eexs}{\end{examples}}

\newcommand{\nul}[1]{\mathrm{N}( {#1} )}
\newcommand{\ran}[1]{\mathrm{R}({#1})}

\newcommand{\norm}[1]{ \| #1 \| }

\newcommand{\codim}{\mathrm{codim \,}}
\newcommand{\dis}[1]{\mathrm{dis }(#1)}

\begin{document}

\title {A new approach in index theory}

\author{ M. Berkani}

\date{}

\maketitle

\begin{abstract}
In this paper, we define an analytical index for a continuous
family of Fredholm operators parameterized by a topological space
$\mathbb{X}$ into a Hilbert space $H,$ as a sequence of integers,
extending naturally the usual definition of the index and  we
prove the homotopy invariance of the index. We give also an
extension of the Weyl theorem for normal continuous families and
we prove that if $H$ is separable, then the space  of B-Fredholm
operators on $H$ is path connected.

\end{abstract}

\footnotetext{\hspace{-7pt}2010 {\em Mathematics Subject
Classification\/}:  47A53, 58B05 \baselineskip=18pt\newline\indent
 }

\section{Introduction}

Let $L(H)$ be  the  Banach algebra  of all bounded  linear
operators defined from an infinite dimensional separable Hilbert
space $H$ to $H,$ $K(H)$ the closed ideal of compact operators on
$H$ and $L(H)/K(H)$ the Calkin Algebra.  We write $\nul{T}$ and
$\ran{T}$ for the  nullspace and the range of an operator $T\in
L(H)$. An operator $T\in L(H)$  is called \cite[Definition
1.1]{BB} a {\it Fredholm} operator if both the nullity of $T,$
$n(T)=\dim \nul{T}$ and the defect of $T, d(T)=\codim \ran{T}$,
are finite. The index $ind(T)$ of a Fredholm operator $T$ is
defined by $ind(T)=n(T)-d(T)$. It is well known that if $T$ is a
Fredholm operator, then $\ran{T}$ is closed.

\bdf \cite{LAB} Let $T\in L(H)$ and let
\[ \Delta (T)=\{ n\in {\mathbb N}: \forall m \in {\mathbb N}, \; m\geq n
\Rightarrow \ran{T^n} \cap \nul{T} \subseteq \ran{T^m} \cap
\nul{T} \}.\] Then the {\it degree of stable iteration of $T$} is
defined as $\dis T= \inf {\Delta (T)}$ (with $\dis {T}=\infty$ if
$\Delta (T)=\emptyset$). \edf

 Define an
equivalence relation $\mathcal{R}$ on the set $fdim(H)\times
fcod(H),$  where $fdim(H)$   is the set of   finite dimensional
vector subspaces of $H$  and $fcod(H)$   is the set of   finite
codimension vector subspaces of $H$  by:

 $$(E_1,F_1) \mathcal{R} (E_1', F_1') \Leftrightarrow  dim E_1 -  codim F_1 = dim E_1'-codim F_1', $$

\noindent Since $H$ is an infinite dimensional vector space, then
the map: $$ \Psi: [fdim(H)\times fcod(H)]/\mathcal{R}  \rightarrow
\mathbb{Z},$$  defined by $\Psi( \overline{(E_1, F_1)} )=dimE_1-
codim F_1,$  where $\overline{(E_1, F_1)} $ is the equivalence
class of the couple $(E_1, F_1),$  is a bijection. Moreover $\Psi$
generate a commutative group structure on the set [$fdim(H)\times
fcod(H)]/R$ and $\psi$ is then a group isomorphism, where $dim$
(resp. $ codim $) is set for the dimension (resp. codimension) of
a vector space.

Reformulating  \cite[Proposition 2.1]{P7}, we obtain.

\bprop \label{prop1} Let $H$ be a Hilbert space and let  $T \in
L(H).$ If there exists an integer $n$ such that  $ (N(T) \cap
R(T^n), R(T) + N(T^n)) $  is  an element of $fdim(H)\times
fcod(H), $ then $d= dis(T)$ is finite and  for all $m \geq d,$  $
(N(T) \cap R(T^m), R(T) + N(T^m)) $ is an element of
$fdim(H)\times fcod(H) $ and
$$ (N(T) \cap R(T^m), R(T) + N(T^m)) \,\,\mathcal{R}\,\, (N(T)
\cap R(T^d), R(T) + N(T^d))$$

\eprop

\bdf \cite[Defintion 2.2]{P7}   Let $H$ be a Hilbert space and let
$T \in L(H).$ Then $T$ is called a B-Fredholm operator if there
exists an integer $n$ such that  $ (N(T) \cap R(T^n), R(T) +
N(T^n)) $  is an element of $fdim(H)\times fcod(H).$  In this case
the index $ ind(T)$ is defined by: $$ ind(T)= dim( N(T) \cap
R(T^n)) - codim(R(T) + N(T^n)) $$  \edf

From Proposition \ref{prop1}, the   definition of the index of a
B-Fredholm  operator is independent of the choice of the integer
$n.$ Moreover,  it extends the usual definition of the index of
Fredholm operators, which are obtained if $n=0.$

 Note also that if  $n$ is an integer  such that  $ (N(T) \cap R(T^n), R(T) +
N(T^n)) $ is an element of $fdim(H)\times fcod(H), $ then from
\cite [Theorem 3.1]{P12}, $R(T^n)$ is closed and the operator $T_n
: R(T^n) \rightarrow R(T^n)$  defined by $T_n(x)= T(x)$ is a
Fredholm operator whose index is equal to the index of $T.$

In  the papers \cite{AS} and \cite{AS1}, the analytical index of a
single or a family of elliptic operators is expressed in terms of
K-theory, and involves in its construction vector bundles.

The first motivation of this work,  is  to give an alternative way
to build an analytical index (called here simply index), for a
continuous family of Fredholm operators parameterized by a
topological space, in the form of a sequence of integers, as a
natural extension of the usual index of a single Fredholm
operator, which is an integer, avoiding the use of vector bundles.

The second motivation  is the extension of  Weyl's theorem
\cite{WL} to continuous families of normal operators.

Moreover, we will prove in the fourth section, that if $H$ is a
separable Hilbert space, then the space $BFred(H)$ is path
connected.

\section{Index of continuous families of Fredholm operators}

\noindent Consider now a family of Fredholm operators parametrized
by a topological space $\mathbb{X},$ that is  a continuous map $
\mathcal{T}: \mathbb{X} \rightarrow Fred(H),$  where $ Fred(H)$ is
the set of Fredholm operators, endowed with the norm topolgy of
$L(H).$ We denote by $ \mathcal{T}_x$ the image $\mathcal{T}(x)$
of an element  $ x \in \mathbb{X}.$

Define   an equivalence relation $ \textbf{\textit{c}} $ on the
space $\mathbb{X}$ by setting that $ x \,\textbf{\textit{c}}\, y,$
if and only if $ x$ and $y$ belongs to the same connected
component of $\mathbb{X}.$ Let $\mathbb{X}^{\textbf{\textit{c}}} $
be the quotient space associated to this equivalence relation and
let $ \mathcal{C}( \mathbb{X}, Fred(H))$ be the space of
continuous maps from the topological space $\mathbb{X}$ into the
topological space $ Fred(H)$ and let  the
 map: $$    q: \mathcal{C}(\mathbb{X}, Fred(H)) \rightarrow [
[fdim(H)\times fcod(H)]/\mathcal{R
}]^{\mathbb{X}^\textbf{\textit{c}}},$$

 defined by $q(\mathcal{T})= (
\overline{N(\mathcal{T}_{x}),R(\mathcal{T}_{x})})_{\overline{x}
\in
 \mathbb{X}^{c}}$ for all  $ \mathcal{T} \in \mathcal{C}( \mathbb{X},
 Fred(H)).$

Define also the map $$\Psi_ {\scriptscriptstyle\mathbb{X}}:
[[fdim(H)\times fcod(H)]/\mathcal{R}
]^{\mathbb{X}^{\textbf{\textit{c}}}}\rightarrow
\mathbb{Z}^{n_\textit{\textbf{c}}}$$ by setting $\Psi_
{\scriptscriptstyle\mathbb{X}}((Y_{\overline{x}}) _{ \overline{x}
\in \mathbb{X}^{\textbf{\textit{c}}} })= (\Psi (Y_{\overline{x}}))
_{\overline{x} \in \mathbb{X}^{\textbf{\textit{c}}}} $ for all $
(Y_{\overline{x}}) _{ \overline{x} \in
\mathbb{X}^{\textbf{\textit{c}}} } \in [fdim(H)\times
fcod(H)]/\mathcal{R} ]^{\mathbb{X}^{\textbf{\textit{c}}}}.$

 Here
$n_{\textbf{\textit{c}}} $ stands for the cardinal of the
connected components of the topological space $\mathbb{X},$
assuming that the space  $\mathbb{X}$ has at most a countable
connected components.

 \bdf \label{index}  The analytical  index (or simply the index) of  a continuous  family of Fredholm operators $ \mathcal{T} : \mathbb{X} \rightarrow Fred(H),$
   parameterized by a topological space $\mathbb{X}$    is defined by $ ind(\mathcal{T})=
\psi_{\scriptscriptstyle  \mathbb{X}}( q(\mathcal{T})).$
 \edf

Explicitly, we have:   $ ind(\mathcal{T})=
\psi_{\scriptscriptstyle {\mathbb{X}}} ((
\overline{N(\mathcal{T}_{x}),R(\mathcal{T}_{x})})_{\overline{x}
\in
 \mathbb{X}^{c}}) = ( ind(\mathcal{T}_x))_{\overline{x}
\in \mathbb{X}^{c}}. $

Thus the index of  a family of B-Fredholm operators $ T$ is a
sequence of integers in $\mathbb{Z}^{n_\textit{\textbf{c}}},$
which may be a finite sequence or infinite sequence, depending on
the cardinal of the connected components of $\mathbb{X}.$

 \bet The index of a continuous family of  Fredholm  operators $\mathcal{T}$  parameterized
by a topological space $\mathbb{X} $ is well defined as an element
of $ \mathbb{Z}^{ n_c}.$ In particular if $\mathbb{X} $ is reduced
to a single element, then the  index of $\mathcal{T}$  is equal to
the usual index of the Fredholm operator $\mathcal{T}.$ \eet

\bp From the usual properties of the index \cite[Theorem
3.11]{BB}, we know that two Fredholm operators located in the same
connected component of the set of Fredholm operators have the same
index. Moreover, as $ \mathcal{T} $ is continuous, the image of a
connected component of the topological space $ \mathbb{X},$  is
included in a connected component of the set of Fredholm
operators.  This shows that the index of a family of Fredholm
operators  is well defined, and it is clear that if $ \mathbb{X}$
is reduced to a single element, the  index of $\mathcal{T}$
defined  here is equal to the usual index of the single Fredholm
operator $\mathcal{T}.$

\ep

\bex Let $ \mathbb{X}= [-1,1] , H= \mathbb{C}$  and $\mathcal{
T}:\mathbb{ X} \rightarrow Fred(\mathbb{C})$ such that
$\mathcal{T}_x(z)= zx $ for all $ z \in \mathbb{C}$ and $ x \in
\mathbb{X}.$ Then from Definition \ref{index}, the index of the
family $\mathcal{T}$ is simply $ind (\mathcal{T})= 0.$ But since
the dimension of the kernel $N(\mathcal{T}_x)$ presents a
discontinuity in $0,$ to build an index using K-theory needs more
tools. \eex

\bdf A continuous  family $\mathcal{K}$  from  $\mathbb{X}$ to
$L(H)$ is said to be compact if $\mathcal{K}_x$ is compact for all
$ x \in \mathbb{X}.$

\edf

\bprop i) Let $\mathcal{T }\in  \mathcal{C}( \mathbb{X}, Fred(H))
$ and  let $ \mathcal{K}$ be a  continuous  compact   family from
$\mathbb{X}$ to $L(H).$ Then $ T+K $ is a Fredholm family and $
ind(T)= ind(T+K).$

ii) Let $\mathcal{S,T}  \in  \mathcal{C}( \mathbb{X}, Fred(H))$ be
two Fredholm families, then the family $\mathcal{ST}$ defined by
$(\mathcal{ST})_x= \mathcal{S}_x\mathcal{T}_x$ is a Fredholm
family and $ ind(\mathcal{ST})= ind(\mathcal{S})+
ind(\mathcal{T}).$

\eprop

\bp This is clear from the usual properties of Fredholm operators.

\ep

\bet \label{Ideal} Assume that $\mathbb{X}$ is a compact
topological space. Then the  set $ \mathcal{K} \mathcal{C}(
\mathbb{X}, L(H))$ of continuous compact families from
$\mathbb{X}$ to $L(H)$ is a closed ideal in the Banach algebra
$\mathcal{C}( \mathbb{X}, L(H))$  \eet

\bp  Recall that $\mathcal{C}( \mathbb{X}, L(H))$ is a unital
algebra with the usual properties of addition, scalar
multiplication and multiplication defined by:

$$ (\lambda \mathcal{S+ T})_x= \lambda \mathcal{S}_x + \mathcal{T}_x,   (\mathcal{ST})_x= \mathcal{S}_x\mathcal{T}_x,  \forall
x \in \mathbb{X}, \lambda \in \mathbb{C}.$$

The unit element of $\mathcal{C}( \mathbb{X}, L(H))$ is the
constant function $\mathcal{I}$ defined by $ \mathcal{I}_x=
I_{\scriptsize H} $ the identity of $H,$  for all $ x \in
\mathbb{X}.$ Moreover as $\mathbb{X}$ is compact, then if we set $
\norm{\mathcal{T}}= sup_{x\in \mathbb{X}} \norm {\mathcal{T}_x},
\forall \mathcal{T} \in \mathcal{C}( \mathbb{X}, L(H)),$  then
$\mathcal{C}( \mathbb{X}, L(H))$ equipped with this norm is a
Banach algebra. Similarly $\mathcal{C}( \mathbb{X}, L(H)/K(H))$
equipped with the norm  $ \norm{\mathcal{T}}= sup_{x\in
\mathbb{X}} \norm {P\mathcal{T}_x}$ is a unital Banach algebra,
where  $ P: L(H) \rightarrow L(H)/K(H)$ is  the usual projection
from $ L(H)$ onto the Calkin algebra $L(H)/K(H).$

 It is clear that $ \mathcal{K} \mathcal{C}( \mathbb{X}, L(H))$ is
 an ideal of  $  \mathcal{C}( \mathbb{X}, L(H)).$ Assume now that
 $ (\mathcal{T}_n)_n$ is a sequence in $ \mathcal{K}\mathcal{C}( \mathbb{X}, L(H))$
 converging  in $ \mathcal{C}( \mathbb{X}, L(H))$
 to $\mathcal{T}.$ Then  $ ((\mathcal{T}_n)_x)_n$ converges to $\mathcal{T}_x,$ as each $(\mathcal{T}_n)_x
 $ is compact, then $\mathcal{T} \in \mathcal{K}\mathcal{C}( \mathbb{X}, L(H)).$
\ep

\bremark  In the same way as in the case of the Calkin algebra,
Theorem \ref{Ideal} generates a new Banach algebra which is
$\mathcal{C}( \mathbb{X}, L(H))/ \mathcal{K} \mathcal{C}(
\mathbb{X}, L(H)).$ Moreover, there is a natural injection
$\overline{\Pi}: \mathcal{C}( \mathbb{X}, L(H))/ \mathcal{K}
\mathcal{C}( \mathbb{X}, L(H)) \rightarrow  \mathcal{C}(
\mathbb{X}, L(H))/K(H))$ defined by $\overline{\Pi}(
\overline{\mathcal{T}})= P\mathcal{T},$ where $
\overline{\mathcal{T}}$ is the equivalence class of the element
$\mathcal{T}$ of $ \mathcal{C}( \mathbb{X}, L(H))$ in
$\mathcal{C}( \mathbb{X}, L(H))/ \mathcal{K} \mathcal{C}(
\mathbb{X}, L(H))$ and $ P: L(H) \rightarrow L(H)/K(H)$ is the
natural projection.

\vspace{3mm} \noindent \textbf{Open question:} Given an element $
\mathcal{S} \in \mathcal{C}(\mathbb{X}, L(H)/K(H)),$  does there
exist a continuous family $ \mathcal{T} \in
\mathcal{C}(\mathbb{X}, L(H))$ such that $ \overline{\Pi}(
\overline{\mathcal{T}})= \mathcal{S}?$

\eremark

\bet
 Assume that $\mathbb{X}$ is a compact topological space  and  let $\mathcal{T} \in \mathcal{C}(\mathbb{X}, L(H)).$ Then $ \mathcal{T}$ is a
 Fredholm family if and only if $P\mathcal{T}$ is  invertible in the Banach algebra
 $ \mathcal{C}(\mathbb{X}, L(H)/K(H)).$

\eet

\bp Assume that $\mathcal{T}$ is a Fredholm family, then for all
$x \in \mathbb{X}, \mathcal{T}_x$  is a Fredholm operator. Thus
$P\mathcal{T}_x$ is invertible in $L(H)/K(H).$ Let
$(P\mathcal{T}_x)^{-1}$ be its inverse, then the family
$(P\mathcal{T})^{-1}$ defined by  $ (P\mathcal{T})^{-1}(x)=
(P\mathcal{T}_x)^{-1}$ is a continuous family, because the
inversion is a continuous map in the Banach algebra $L(H)/K(H),$
and $(P\mathcal{T})^{-1}$ is the inverse of $P\mathcal{T}$ in the
Banach algebra $ \mathcal{C}(\mathbb{X}, L(H)/K(H)).$

Conversely if $ P\mathcal{T}$  is  invertible in the Banach
algebra
 $ \mathcal{C}(\mathbb{X}, L(H)/K(H)),$ then there exists $\mathcal{S} \in   \mathcal{C}(\mathbb{X}, L(H)/K(H))$
such that $ (P\mathcal{T})\mathcal{S}= \mathcal{S}(P\mathcal{T})=
\overline{\mathbb{\mathcal{I}}},$  where
$\overline{\mathbb{\mathcal{I}}}$ is defined by
$\overline{\mathbb{\mathcal{I}}}_x= PI_H,$ for all $ x \in
\mathbb{X}, \,I_H$ being the identity of $H.$ Thus
$(P\mathcal{T}_x) \mathcal{S}_x= \mathcal{S}_x(P\mathcal{T}_x)=
PI_H.$ Thus $P\mathcal{T}_x $ is invertible in the Calkin algebra
$L(H)/K(H),$
 $\mathcal{T}_x$ is a Fredholm operator and $\mathcal{T} \in \mathcal{C}(\mathbb{X},
Fred(H)).$

\ep


\bdf Let $ \mathcal{S,T} $   be in $\in \mathcal{C}( \mathbb{X},
Fred(H)). $ We will say that $\mathcal{S}$ and $\mathcal{T}$ are
Fredholm homotopic, if there exists a map $\Phi:  [0,1]\times X
\rightarrow L(H)$ such for all $ (t,x) \in [0,1]\times
\mathbb{X},\, \Phi(0,x)= \mathcal{S}_x, \, \Phi(1,x)=
\mathcal{T}_x $ and $\Phi(t,x)$ is a Fredholm operator. \edf

\bet Let $ \mathcal{S, T}$ be two Fredholm homotopic elements of
\,$\mathcal{C}( \mathbb{X}, Fred(H)). $  Then  $ind(\mathcal{T})=
ind(\mathcal{S}).$ \eet

\bp Since $\mathcal{S}$ and $\mathcal{T}$ are Fredholm homotopic,
there exists a continuous map $ h: \mathbb{X} \times [0, 1]
\rightarrow Fred(H) $  such that such that $h(x,0)=\mathcal{S}(x)$
and $h(x,1)=\mathcal{T}(x)$ for all $x\in\mathbb{X}.$ For a fixed
$x \in \mathbb{X} ,$ the map $ h_x: [0,1] \rightarrow Fred(H),$
defined by $h_x(t)= h(x,t)$ is a continuous path
 in $ Fred(H)$ linking  $\mathcal{S}_x$ to $\mathcal{T}_x.$  Thus $ind(\mathcal{S}_x)=
 ind(\mathcal{T}_x).$ So  $q(\mathcal{S})=q(\mathcal{T})$ and then  $ind(\mathcal{S})= ind(\mathcal{T}).$

\ep

\bet Let $\mathbb{X}$ be a compact  topological  space. Then the
index is a continuous locally constant function  from $
\mathcal{C}( \mathbb{X}, Fred(H)) $  into the group $\mathbb{Z}^{
n_\textit{\textbf{c}}}.$

\eet

\bp Let $\mathcal{T} \in \mathcal{C}( \mathbb{X}, Fred(H)),$  then
$\forall x \in  \mathbb{X}, \exists \, \epsilon_x > 0, $ such that
$B(\mathcal{T}_x, \epsilon_x) \subset Fred(H),$  because $Fred(H)$
is open in $L(H).$ Then the index is constant on $B(\mathcal{T}_x,
\epsilon_x),$ because $B(\mathcal{T}_x, \epsilon_x)$ is connected.
We have $ \mathbb{X} \subset \bigcup \limits_{x\in \mathbb{X}}
\mathcal{T}^{-1}( B(\mathcal{T}_x, \frac{\epsilon_{x}}{2})).$
Since $\mathbb{X}$ is compact, there exists $ x_1,...x_n$ in
$\mathbb{X}$ such that $\mathbb{X} \subset \bigcup \limits_{i=1}^n
\mathcal{T}^{-1}(B(\mathcal{T}_{x_i}, \frac{\epsilon_{x_i}}{2})).$
Let $ \epsilon= min \{ \frac{\epsilon_{x_i}}{2} | 1\leq i \leq n
\} $ the minimum of the $ \frac{\epsilon_{x_i}}{2}, 1\leq i \leq
n,$ and  let $\mathcal{S} \in \mathcal{C}( \mathbb{X}, Fred(H)) $
such that $ || \mathcal{T-S}|| < \frac{\epsilon}{2}.$ If $ x \in
\mathbb{X},$ then $ ||\mathcal{ T}_x-\mathcal{S}_x|| <
\frac{\epsilon}{2}$ and there exists $ i, 1\leq i \leq n,$ such
that $ x \in \mathcal{T}^{-1}(B(\mathcal{T}_{x_i},
\frac{\epsilon_{x_i}}{2})).$ Then $ || \mathcal{S}_x
-\mathcal{T}_{x_i}|| \leq || \mathcal{S}_x-\mathcal{T}_x||  +
||\mathcal{T}_x- \mathcal{T}_{x_i}|| < \epsilon/2 +
\epsilon_{x_i}/2 \leq \epsilon_{x_i}.$ So $ind(\mathcal{S}_x)=
ind(\mathcal{T}_{x_i})= ind(\mathcal{T}_x).$ Hence the index is a
locally constant function, in particular it is a continuous
function. \ep

\bet Let $\mathbb{X}$ be a compact topological space.  Then the
set $ \mathcal{C}( \mathbb{X}, Fred(H))$ is an open subset of the
Banach algebra $ \mathcal{C}( \mathbb{X}, L(H))$ endowed with the
uniform norm $ \norm{ \mathcal{T}}= sup_{x\in  \mathbb{X}}
\norm{\mathcal{T}_x}.$ \eet

\bp Let $\mathcal{T} \in \mathcal{C}( \mathbb{X}, Fred(H)),$  then
$\forall x \in \mathbb{X}, \exists \, \epsilon_x
> 0, $ such that $B(\mathcal{T}_x, \epsilon_x) \subset Fred(H).$   We have $ \mathbb{X} \subset
\bigcup\limits_{x\in \mathbb{X}} \mathcal{T}^{-1}(
B(\mathcal{T}_x, \frac{\epsilon_{x}}{2})).$ Since $\mathbb{X}$ is
compact, there exists $ x_1,...x_n$ in $\mathbb{X}$ such that
$\mathbb{X} \subset \bigcup \limits_{i=1}^n
\mathcal{T}^{-1}(B(\mathcal{T}_{x_i}, \frac{\epsilon_{x_i}}{2})).$
Let $ \epsilon= min \{ \frac{\epsilon_{x_i}}{2} | 1\leq i \leq n
\} $ the minimum of the $ \frac{\epsilon_{x_i}}{2}, 1\leq i \leq
n.$ Let $\mathcal{S} \in \mathcal{C}( \mathbb{X}, L(H)) $ such
that $ || \mathcal{T-S}|| < \frac{\epsilon}{2}.$ If $ x \in
\mathbb{X},$ then $ || \mathcal{T}_x-\mathcal{S}_x|| <
\frac{\epsilon}{2}$ and there exists $ i, 1\leq i \leq n,$ such
that $ x \in \mathcal{T}^{-1}(B(\mathcal{T}_{x_i},
\frac{\epsilon_{x_i}}{2})).$ Then  $ || \mathcal{S}_x
-\mathcal{T}_{x_i}|| \leq || \mathcal{S}_x-\mathcal{T}_x|| +
||\mathcal{T}_x- \mathcal{T}_{x_i}|| < \epsilon/2 +
\epsilon_{x_i}/2 \leq \epsilon_{x_i}$ and $\mathcal{S}_x$ is a
Fredholm operator. Therefore $\mathcal{S} \in \mathcal{C}(
\mathbb{X}, Fred(H))$ and $\mathcal{C}( \mathbb{X}, Fred(H))$ is
open in $ \mathcal{C}( \mathbb{X}, L(H)).$

Alternatively, we can see that $ \mathcal{C}( \mathbb{X},
Fred(H))= \Pi^{-1}( (\mathcal{C}( \mathbb{X}, L(H)/K(H)) ^{inv}),$
where $(\mathcal{C}( \mathbb{X}, L(H)/K(H)) ^{inv}$ is the open
group of invertible elements of the unital Banach algebra
$\mathcal{C}( \mathbb{X}, L(H)/K(H))$ and $ \Pi: \mathcal{C}(
\mathbb{X}, L(H))\rightarrow \mathcal{C} ( \mathbb{X}, L(H)/K(H))$
is the map defined by $ \Pi(\mathcal{T})= P\mathcal{T},$ for all $
\mathcal{T} \in \mathcal{C}( \mathbb{X}, L(H)/K(H)).$

\ep

\vspace{2mm}

If $T$ is a Fredholm operator of index $0,$  then $T$ is Fredholm
connected to the identity $ I_H$ of $H.$ However,  it follows from
\cite[Remark 3.32]{BB}, that there exists continuous Fredholm
families whith vanishing index, which are not Fredholm homotopic
to the identity $\mathcal{I}$  of $ \mathcal{C}( \mathbb{X},
L(H)).$  This is the reason why we introduce the following class
of topological spaces.

\bdf Let $\mathbb{X}$ be a  compact topological space. We will say
that $\mathbb{X}$ satisfies the $\mathcal{H}$-condition, if every
continuous Fredholm family of index $0$ in $\mathbb{Z}^{n_c},$  is
Fredholm homotopic to the identity  $\mathcal{I}$  of $
\mathcal{C}( \mathbb{X}, L(H)).$  \edf

\bexs \begin{enumerate}

\item Let $\mathbb{X}$ be a compact contractible space. Then there
exists $x_0 \in X,$  a continuous map $\Phi: [0,1]\times
\mathbb{X} \rightarrow \mathbb{X},$ such that $\Phi(x,0)= x,
\Phi(x,1)= x_0,$ for all $x \in \mathbb{X}.$

If  $\mathcal{T}$ in $\mathcal{C}(\mathbb{X}, L(H))$ is a
continuous Fredholm family of  index $0,$ then $\mathcal{T}_{x_0}
$ is a Fredholm operator of index $0.$ So $\mathcal{T}_{x_0} $ is
Fredholm connected to the identity $ I_H$ of $H.$ Hence  there
exist a map $\varphi: [0, 1] \rightarrow Fred(H),$ such that
$\varphi(0)= \mathcal{T}_{x_0},$ and $\varphi(1)= I_H.$

Define the  map $\Psi: X\times [0, 1] \rightarrow Fred(H),$ by
setting  $\Psi( x,t)= \varphi(t),$ for all $(x,t) \in
\mathbb{X}\times [0,1].$ Since $\varphi$ is continuous, $\Psi$ is
continuous.

Now let  $\Gamma : X\times [0, 1] \rightarrow Fred(H),$ defined by
setting $\Gamma( x,t)=\mathcal{ T} (\Phi(x,2t)), $ for all $(x,t)
\in \mathbb{X}\times [0,\frac{1}{2}],$ and  $\Gamma( x,t)=
\Psi(x,2t-1),$ for all $(x,t) \in \mathbb{X}\times
[\frac{1}{2},1].$ Then  $\Gamma$ is a continuous map from $X\times
[0, 1] \rightarrow Fred(H),$ satisfying
$\Gamma(x,0)=\mathcal{T}_x$ and $\Gamma(x,1)=I_H,$ for all $x \in
\mathbb{X}.$

Hence $ \mathcal{T}$ is Fredholm homotopic to the identity
$\mathcal{I}$  and $\mathbb{X} $ satisfies the
$\mathcal{H}$-condition.

 \item A disjoint union of
compact convex sets  in a topological vector space,  satisfies the
$\mathcal{H}$-condition, as each convex set is contractible.
\end{enumerate}

\eexs

The $\mathcal{H}-$condition is a necessary condition in order to
have that the  index is a bijective monoid isomorphism map from
the homotopy equivalence classes $[\mathbb{X}, Fred(H)]$ of
continuous families of Fredholm operators into the group
$\mathbb{Z}^{ n_\textit{\textbf{c}}},$ similarly to the
Atiyah-J$\ddot{a}$nich  theorem \cite[Theorem 3.40]{BB}. In the
case of a locally connected compact topological space, we show
that the $\mathcal{H}$-condition is   also sufficient.

\bet \label{iso} Let $\mathbb{X}$ be a locally connected compact
topological space. Then the index is a bijective monoid
isomorphism map from the homotopy equivalence classes
$[\mathbb{X}, Fred(H)]$ of continuous families of Fredholm
operators into the group $\mathbb{Z}^{ n_\textit{\textbf{c}}}$ if
and only if  $\mathbb{X}$ satisfies the $\mathcal{H}$-condition.

\eet

\bp  Observe first that the index is a monoid homomorphism, since
$ ind(\mathcal{I})= 0,$ and $ ind (\mathcal{S}\mathcal{T})=
ind(\mathcal{S})+ ind(\mathcal{T}),$  for all $\mathcal{S},
\mathcal{T}$ in  $\mathcal{C}(\mathbb{X}, Fred(H)).$

Assume  now that  the index is a bijective monoid isomorphism map
from the homotopy equivalence classes $[\mathbb{X}, Fred(H)]$ of
continuous families of Fredholm operators into the group
$\mathbb{Z}^{ n_\textit{\textbf{c}}}.$ If $\mathcal{T}$ is a
continuous Fredholm family of index $0, $ then $\mathcal{T}$ is
Fredholm homotopic to the identity $\mathcal{I},$ because the
indentity $\mathcal{I}$ is of index $0.$

Conversely   assume that $\mathbb{X}$ satisfies the
$\mathcal{H}$-condition and let $\mathcal{T} \in \mathcal{C}(
\mathbb{X}, Fred(H)) $  such that $ ind( \mathcal{T}) = 0. $ Then
the Fredholm family $\mathcal{T}$ is Fredholm homotopic to the
indentity $\mathcal{I}$.  Thus the index is injective.

As   $\mathbb{X}$ is a  locally connected compact topological
space, its connected components are clopen subsets of $\mathbb{X}$
and their cardinal ${n_c}$ is finite. Let $ C_{1}, ... C_{n_c} $
be those connected components.  To prove that the index is
surjective, let $u = (n_1,n_2,...,n_{n_c}) \in \mathbb{Z}^{n_c}.$
On each connected component of $\mathbb{X}$, take a Fredholm
operator $T_{i}$ such that $ ind(T_{i})= n_i$, for $ 1\leq i \leq
n_c $ and define the map $ \mathcal{T}: \mathbb{X} \rightarrow
\mathbb{Z}^ {n_c}$  by $ \mathcal{T}(x)= T_i,$ if $ x \in C_{i}.$
Then it is clear that $\mathcal{T}$ is continuous and that $
ind(\mathcal{T})= (n_1,n_2,...,n_{n_c}).$ Thus the index is
surjective.

 \ep

\vspace{2mm}

 We can conclude from Theorem \ref{iso}, that in the
case of a locally connected, compact topological space satisfying
the $\mathcal{H}$-condition, the index defined here, is up to an
isomorphism the same as the index bundle defined in \cite{AS} and
\cite{AS1}.

\section{ An application to Weyl theorem}
In this section, by an application to Weyl's theorem,  we show
that the index of continuous Fredholm families defined in section1
is a natural generalization of the usual index of a single
Fredholm operator.

So we will extend  Weyl theorem to continuous normal families of
bounded linear operators.  Recall that the classical Weyl's
theorem \cite{WL} asserts that if $T$ is a normal operator acting
on a Hilbert space $H$, then the Weyl spectrum $\sigma_W(T)$ is
exactly the set of all points in $\sigma(T)$ except the isolated
eigenvalues of finite multiplicity, that is
 $$ \sigma_W(T)= \sigma(T) \backslash
E_0(T),$$  where $E_0(T)$ is the set of isolated eigenvalues of
$T$ of  finite multiplicity and $\sigma_W(T)$ is the Weyl spectrum
of $T$, that is $\sigma_{W}(T)=\{ \lambda \in {\bf C} \mid  T -
\lambda I $ is not a Fredholm operator of index $0$\}.

\bdf  Let $\mathcal{T} \in \mathcal{C}( \mathbb{X}, L(H)). $ Then
$\mathcal{T}$ is called a Weyl   family
 if it is a
Fredholm   family of index $0$ in $\mathbb{Z}^{n_c}.$\\
  \noindent The  Weyl
spectrum $\sigma_{W}(\mathcal{T})$ of $\mathcal{T}$ is defined by
$\sigma_{W}(\mathcal{T})= \{ \lambda \in {\bf C}: \mathcal{T}-
\lambda \mathcal{I} $ is not a Weyl family \}.

\edf

\bdf Let $T \in L(H).$

\begin{enumerate}

\item The  ascent  $a(T)$ of  $T$ is defined by
 $a(T)=\mbox{inf} \{ n\in \mathbb{N} \mid N(T^n)=N(T^{n+1})\},$
and the  descent   $ \delta(T)$ of $T$, is defined by
 $\delta(T)= \mbox{inf} \{ n \in \mathbb{N} \mid R(T^n)= R(T^{n+1})\},$ with $ \mbox{inf}\, \emptyset= \infty.$

 \item  A complex number $ \lambda$ is a pole of the
resolvent of $T$ of finite rank  if  $ 0 < n(T-\lambda I) < \infty
\,\,$    and max $(a(T- \lambda I), \delta(T- \lambda I))<
\infty.$ Moreover, if this is true, then $a(T- \lambda I)=
\delta(T- \lambda I).$

\end{enumerate}

\edf

\bdf  Let $\mathcal{T} \in \mathcal{C}( \mathbb{X}, L(H)),$ then:

\begin{enumerate}

\item $E_0(\mathcal{T})$ is defined by  $  E_0( \mathcal{T})= \{
\lambda \in \sigma(\mathcal{T}) \mid \exists A\subset \mathbb{X},
A \neq \emptyset: \lambda \in \bigcap \limits _{x \in
\mathbb{A}}{} E_0 ( \mathcal{T}_x) , \lambda \notin \bigcup
\limits _{x \in \mathbb{\mathbb{X}\setminus A}}{}\sigma(
\mathcal{T}_x) \},$
 where $
E_0(\mathcal{\mathcal{T}}_x)$ is the set of isolated eigenvalues
of finite multiplicity of  $ {\mathcal{T}}_x $  in
$\sigma{(\mathcal{T}}_x).$

\item $\Pi_0(\mathcal{T})$ is defined by  $ \Pi_0( \mathcal{T})=
\{  \lambda \in \sigma(\mathcal{T}) \mid \exists A\subset
\mathbb{X}, A \neq \emptyset: \lambda \in \bigcap \limits _{x \in
\mathbb{A}}{} \Pi_0 ( \mathcal{T}_x) , \lambda \notin \bigcup
\limits _{x \in \mathbb{\mathbb{X}\setminus A}}{}\sigma(
\mathcal{T}_x) \},$ where $ \Pi_0(\mathcal{\mathcal{T}}_x)$ is the
set of poles of $ {\mathcal{T}}_x$  of finite rank.

\end{enumerate}

\edf

\bdf   \label{def:weyl-theo} Let $\mathcal{T} \in \mathcal{C}(
\mathbb{X}, L(H)).$  We will say that:

\begin{enumerate}

\item $\mathcal{T}$ satisfies  Browder's theorem  \/  if
$\sigma_{W}(\mathcal{T}) = \sigma(\mathcal{T}) \setminus
\Pi_{0}(\mathcal{T}).$

\item $\mathcal{T}$ satisfies Weyl's theorem \/   if
$\sigma_{W}(\mathcal{T}) = \sigma(\mathcal{T}) \setminus
E_{0}(\mathcal{T}).$

\end{enumerate}

\edf

\bet \label{BW}

 Let $\mathcal{T} \in \mathcal{C}( \mathbb{X}, L(H)).$ Then the
 following holds:

\begin{enumerate}

\item If $\mathcal{T}_x$ satisfies Weyl's theorem for all $ x \in
\mathbb{X},$ then    $\mathcal{T}$ satisfies Weyl's theorem.

\item If $\mathcal{T}_x$ satisfies Browder's  theorem for all $ x
\in \mathbb{X},$   then $\mathcal{T}$ satisfies Browder's theorem.

\end{enumerate}

\eet

\bp 1- First, let us show that  $\sigma_{W}(\mathcal{T})=  \bigcup
\limits _{x \in \mathbb{X}}{} \sigma_{W}(\mathcal{T}_x).$

\noindent If $\lambda \notin \sigma_{W}(\mathcal{T}),$ then
$\forall x \in \mathbb{X}, ind(\mathcal{T}_x - \lambda I) =0.$
Thus $ \lambda \notin \bigcup \limits _{x \in \mathbb{X}}{}
\sigma_{W}(\mathcal{T}_x).$

\noindent Conversely if $ \lambda \notin \bigcup \limits _{x \in
\mathbb{X}}{} \sigma_{W}(\mathcal{T}_x),$  then  $\forall x \in
\mathbb{X}, ind(\mathcal{T}_x - \lambda I) =0.$ So on each
connected component $\mathcal{C}$ of $ \mathbb{X}$, we have
$ind(\mathcal{T}_x - \lambda I)=0 , \forall x \in \mathcal{C}.$
 Using the definition of the index of a Fredholm
family, we obtain that $ind(\mathcal{T}_x - \lambda I) =0$ in
$\mathbb{Z}^{n_c}.$

\noindent Moreover as the inversion is a continuous map on the
group of invertible elements of $L(H),$ then $\sigma(\mathcal{T})=
\bigcup \limits _{x \in \mathbb{X}}{} \sigma(\mathcal{T}_x).$

Assume now that $ \lambda \in \sigma(\mathcal{T})$ and  $ \lambda
\notin \sigma_W(\mathcal{T}).$ Then $\forall x \in \mathbb{X},
\lambda \notin \sigma_W(\mathcal{T}_x).$   If $ \lambda \notin
E_0(\mathcal{T}_x),$  as $ \mathcal{T}_x$ satisfies Weyl's
theorem, then   $\lambda \notin \sigma(\mathcal{T}_x).$ As $
\lambda \in \sigma(\mathcal{T}),$ then there exists $ x_0 \in
\mathbb{X}$ such that  $ \lambda \in \sigma(\mathcal{T}_{x_0}).$
So $ \lambda \in  E_0(\mathcal{T}_{x_0}).$
 Let $A =  \{ x \in \mathbb{X}
\mid \lambda \in  E_0 ( \mathcal{T}_x)\},$  then $ A \neq
\emptyset.$ Moreover  $ \lambda \in \bigcap \limits _{x \in
\mathbb{A}}{} E_0 ( \mathcal{T}_x)$ and $ \lambda \notin \bigcup
\limits _{x \in \mathbb{\mathbb{X}\setminus A}}{}\sigma(
\mathcal{T}_x).$ So $ \lambda \in E_0(\mathcal{T})$  and
$\sigma(\mathcal{T}) \subset  \sigma_{W}(\mathcal{T}) \bigcup
E_0(\mathcal{T}).$ As we have always  $\sigma_{W}(\mathcal{T})
\bigcup E_0(\mathcal{T})  \subset  \sigma(\mathcal{T}), $ then
$$\sigma(\mathcal{T}) = \sigma_{W}(\mathcal{T}) \bigcup
E_0(\mathcal{T}).$$

Suppose now that $ \lambda \in \sigma_{W}(\mathcal{T}) \bigcap
E_0(\mathcal{T}).$ Then   $A =  \{ x \in \mathbb{X} \mid \lambda
\in E_0 ( \mathcal{T}_x)\}\neq \emptyset$ and $\forall x \in
A,\lambda \notin \sigma_{W}(\mathcal{T}_x),$ because $
\mathcal{T}_x$ satisfies Weyl's theorem.  If $x \in \mathbb{X}
\setminus A,$ then $\lambda \notin \sigma(\mathcal{T}_x)$ and then
$\lambda \notin \sigma_W(\mathcal{T}_x).$ So $ \lambda \notin
  \bigcup \limits _{x \in \mathbb{X}}{}
\sigma_{W}(\mathcal{T}_x)= \sigma_{W}(\mathcal{T}).$ But this is a
contradiction because $ \lambda \in \sigma_{W}(\mathcal{T}).$
Therefore  $\sigma_{W}(\mathcal{T}) \bigcap E_0(\mathcal{T})=
\emptyset.$ Finally we have  $\sigma_{W}(\mathcal{T}) =
\sigma(\mathcal{T}) \setminus E_{0}(\mathcal{T})$ and so  $
\mathcal{T}$ satisfies Weyl's theorem

2- We use the same proof as in the first part, just by  replacing
$ E_0( \mathcal{T})$ by  $\Pi_0( \mathcal{T}).$

\ep

\bdf  Let $\mathcal{T} \in \mathcal{C}( \mathbb{X}, L(H)).$  Then
$\mathcal{T}$ is called a  normal family if $ \forall x \in
\mathbb{X}, \mathcal{T}_x $ is a normal operator. \edf

We give now an extension  of  Weyl's theorem to the case of normal
families.

\bet  Let $\mathcal{T} \in \mathcal{C}( \mathbb{X}, L(H))$ be a
normal family. Then $\mathcal{T}$  satisfies Weyl's theorem,
that's $\sigma_{W}(\mathcal{T})= \sigma( \mathcal{T}) \setminus
E_0(\mathcal{T}).$ \eet

\bp Since $\mathcal{T}$ is a normal family, then $ \forall x \in
\mathbb{X}, \mathcal{T}_x $ is a normal operator. From \cite{WL},
we know that a normal operator satisfies Weyl's theorem.
 The theorem is then a consequence of Theorem \ref{BW}. \ep

\section{The  space of  B-Fredholm operators}

In this section, we  prove that the space of B-Fredholm operators
$BFred(H)$  on a separable Hilbert space  $H$  is path connected,
and  the index function is not continous on  $BFred(H)$ .

\bdf Let $ S,T $   be in $  BFred(H). $ We will say that $S$ and
$T$ are B-Fredholm homotopic, if there exists a continuous map
$\Phi: [0,1]\rightarrow BFred(H)$  such  $ \Phi(0)= S $ and
$\Phi(1)= T.$ \edf

\bprop \label{B-Fredh-connected} Let $S, T$ be two B-Fredholm
operators acting on a Hilbert space having equal index. Then $S$
and $T$ are B-Fredholm path connected.

\eprop

\bp Let $S,  T \in BFred(H)$ such that $ind(S)= ind(T)= n.$ Then
From \cite[Remark A]{P12} there exists $\epsilon
> 0$ such if $\lambda \in \mathbb{C}, $   $ 0<\mid \lambda\mid <
\epsilon,$ both of  $S- \lambda I$ and $T- \lambda I$ are Fredholm
operator and $ind( S- \lambda I)= ind( T- \lambda I)= n.$

Let  $\lambda$ such that $0<\mid \lambda\mid < \epsilon.$   Since
$ S- \lambda I$ and $ T- \lambda I $  are Fredholm operators
having the same index, then they are Fredholm path-connected.

Consider now the map $\phi: [0, 1] \rightarrow L(H),$ defined by $
\phi(t)= S- t\lambda I.$ Then $ \phi$ is continuous, $ \phi(0)= S,
\phi(1)= S - \lambda I $ and for all $ t \in [ 0, 1], \phi(t)$ is
a B-Fredholm operator. Thus $ S $ is  B-Fredholm path-connected to
$S- \lambda I.$ Similarly $T $ is B-Fredholm path connected to $T-
\lambda I.$ By the transitivity of the path-connectedness, $S$ and
$T$ are B-Fredholm path-connected.

\ep

We prove now  the main result of this section.

\bet \label{T-BP}

Let $H$ be a separable Hilbert space. Then the topological space
$BFred(H)$ is path connected. \eet

\bp Without loss of generality, we may assume  that
$H=l^2(\mathbb{C}).$  Let $S, T$ be the operators defined on $H$
by:

$ S(x_1,x_2,...,x_n,...)=(x_1,0,0,0,..,0,..), \forall x=(x_i)_i
\in l^2(\mathbb{C}),$

$T(x_1,x_2,...,x_n,..)=(x_1,x_3,x_4,x_5,x_6,...), \forall
x=(x_i)_i \in l^2(\mathbb{C}). $

Then $S$ is a B-Fredholm operator of index $0$ and $T$ is a
Fredholm operator of index $1.$ Let $ \Phi: [0,1]  \rightarrow
L(H)$ be the map defined by:

 $\Phi(t)(x_1,x_2,...,x_n,...)= ( x_1, tx_3,
tx_4,tx_5,....),{\text for \,  all}\,\, x=(x_i)_i \in
l_2(\mathbb{C}).$

Then $ \Phi$ is continuous and if $ 0<t\leq 1,$
 $\Phi(t)$ is a Fredholm operator of index $1.$ So $ \Phi$ is a continuous
 path of B-Fredholm operators such that
$ \Phi(0)= S $ and $\Phi(1)= T.$  Thus $S$ and $T$ are B-Fredholm
path  connected.

As $S$ is of index $0$, then from Proposition
\ref{B-Fredh-connected}, $S$ is B-Fredholm path connected to the
identity, because they have the same index. Thus $T$ is B-Fredholm
path connected to the identity operator $I$. Consequently, $
\forall n \in \mathbb{N}, T^n$ is path connected to the identity
operator $I= I^n.$ We observe that $ ind(T^n)= n ind(T)=n.$ Using
the adjoint, we deduce that $ {T^*}^n$ is B-Fredholm path
connected to the identity operator $I,$  with $ ind({T^*})=-1$ and
$ ind({T^*}^n)= -n.$

Let $U, V$ be two B-Fredholm operators of indexes $m \geq 0$ and
$p\leq 0$ respectively. Again from Proposition
\ref{B-Fredh-connected}, $U$ is B-Fredholm path connected to $T^m$
and $V$ is B-Fredholm path connected  to $(T^*)^{-p}.$ Hence $U$
and $V$ are B-Fredholm path connected to $I,$  and so they are
B-Fredholm path connected. From here, we conclude that any two
B-Fredholm operators are B-Fredholm path connected and  the space
$BFred(H)$ is path connected. But the index function is highly
discontinuous on $BFred(H).$

\ep

As a consequence of  Theorem \ref{T-BP}, we cannot define an
analytical index for continuous B-Fredholm families as done for
continuous Fredholm families.

\vspace{3mm}

{  \bf  Acknowledgment:} The author would like to thank the
referee for his comments, which led to an important improvement of
Theorem\ref{iso}.

 \baselineskip=12pt
\bigskip
\vspace{-5 mm }
 \baselineskip=12pt
\bigskip

{\small
\noindent Mohammed Berkani,\\
 \noindent Science faculty of Oujda,\\
\noindent University Mohammed I,\\
\noindent Laboratory LAGA, \\
\noindent Morocco\\
\noindent berkanimo@aim.com\\

\end{document}